\documentstyle[11pt]{article}

\begin{document}

\title{Revised support forcing iterations and the continuum hypothesis}
\author{{\it Preliminary report} \\
Chaz Schlindwein \\
CSC 690, Spring 1996}

\maketitle

\def\forces{\mathbin{\parallel\mkern-9mu-}}
\def\notforces{\,\nobreak\not\nobreak\!\nobreak\forces}
\def\Gen{{\rm Gen}}
\def\calm{{\cal M}}
\def\bfone{{\bf 1}}

\def\restr{\,\hbox{\vrule height8pt width.4pt depth0pt
\vrule height7.75pt width0.3pt depth-7.5pt\hskip-.2pt
\vrule height7.5pt width0.3pt depth-7.25pt\hskip-.2pt
\vrule height7.25pt width0.3pt depth-7pt\hskip-.2pt
\vrule height7pt width0.3pt depth-6.75pt\hskip-.2pt
\vrule height6.75pt width0.3pt depth-6.5pt\hskip-.2pt
\vrule height6.5pt width0.3pt depth-6.25pt\hskip-.2pt
\vrule height6.25pt width0.3pt depth-6pt\hskip-.2pt
\vrule height6pt width0.3pt depth-5.75pt\hskip-.2pt
\vrule height5.75pt width0.3pt depth-5.5pt\hskip-.2pt
\vrule height5.5pt width0.3pt depth-5.25pt}\,}

\def\overtau{{\overline\tau}}

\noindent Current address (until June 15):

\noindent 995 Count Wutzke Av

\noindent Las Vegas NV 89119  USA

\noindent {\tt chaz@cs.unlv.edu}

\medskip

\noindent After June 15:

\noindent Division of Mathematics and Computer Science

\noindent Lander University

\noindent Greenwood SC 29649 USA

\medskip

\centerline{{\bf Abstract}}

\def\cf{{\rm cf}}
\def\sup{{\rm sup}}
\def\supt{{\rm supt}}
\def\dom{{\rm dom}}
\def\range{{\rm range}}
\def\calj{{\cal J}}

In [Sh, chapter XI], Shelah shows that certain revised
countable support (RCS) iterations
do not add reals.  His motivation is to establish the
independence (relative to large cardinals) of Avraham's
problem on the existence of uncountable non-constuctible
sequences all of whose proper initial segments are constructible,
Friedman's problem on whether every 2-coloring of
$S^2_0=\{\alpha<\omega_2\colon\cf(\alpha)=\omega\}$
has an uncountable sequentially closed homogeneous
subset, and existence of a precipitous normal filter on $\omega_2$
with $S^2_0\in{\cal F}$.
The posets which Shelah uses in these constructions are Prikry forcing,
Namba forcing, and the forcing consisting of closed
countable subsets of $S^*$ under reverse end-extension,
where $S^*$ is a fixed stationary co-stationary subset of $S^2_0$.
Shelah establishes different preservation theorems for each of these
three posets (the theorem for Namba forcing is
particularly intricate).  We establish a general
preservation theorem for a variant of RCS iterations which
includes all three posets in a straightforward way.

\vfill\eject

Revised countable support (RCS) iterations, introduced by Shelah
{Sh, chapter X], is a generalization of the well-known countable
support iteration method for forming forcing iterations.
RCS iterations are appropriate for semi-proper forcings (e.g.,
Prikry forcing) and also even for forcing iterations built of
posets which cofinalize some regular cardinal to have countable
cofinality even if the forcings are not semi-proper (for example,
one may use RCS iteration to iterate Namba forcing even in cases
where Namba forcing is not semi-proper). The most important
fact about RCS iterations is that semi-properness is preserved
(i.e., if the constituent posets of the iteration are semi-proper,
then the resulting forcing is itself semi-proper).
This was established by Shelah [Sh, chapter X]. A simplified version
of RCS iterations was introduced in [S] along with a simplified
proof of preservation of semi-properness.

A secondary fact about RCS iterations is that certain posets can
be iterated without adding reals.  For example, roughly speaking,
Shelah [Sh, chapter XI] shows that
iteration of Prikry forcing, or iteration of Namba forcing over a
ground model satisfying CH, does not add reals (regardless of
whether Namba forcing is semi-proper).  Certain technical restrictions
apply in 
Shelah's results; notably, in the case of 
Namba forcing certain cardinal collapses must occur, via
Levy collapses or other $\sigma$-closed
cardinal collapses.  
In this paper we introduce a variant of revised countable
support and prove a general preservation theorem for not
adding reals in which many of
the technical restrictions of Shelah are eliminated.

{\bf Notation.}  Our notation follows [S]. We set
$M^P$ equal to the set (or class) of $P$-names which are in $M$.
This is different from the class of names whose values are
forced to be in $M$, and it is different from the class of
names whose values are forced to be in $M[G_P]$ (of course, for
any $G$ and any name $\dot x$ which is forced to be in $M[G_P]$
there is $p\in G$ and $\dot y \in M^P$ such that $p\forces``\dot x
=\dot y$'').
The notation $\dot P_{\eta,\alpha}$ is used in the context
of a forcing iteration $\langle P_\beta\colon\beta\leq\kappa\rangle$
based on $\langle\dot Q_\beta\colon\beta<\kappa\rangle$; it
denotes a $P_\eta$-name characterized by
$V[G_{P_\eta}]\models``\dot P_{\eta,\alpha}=\{p\restr[\eta,\alpha)
\colon p\restr\eta\in G_{P_\eta}$ and $p\in P_\alpha\}$.''
Here, $s\restr[\eta,\alpha)$ is not the check (with respect to
$P_\eta$) of the restriction of $s$ to the interval
$[\eta,\alpha)$, but rather it is a $P_\eta$-name for the function
$f$ with domain equal to $[\check\eta,\check\alpha)$
such that $f(\check\beta)$ is the $P_\eta$-name for the
$\dot P_{\eta,\beta}$-name corresponding to the $P_\beta$-name
$s(\beta)$.  This may be contrasted with the definition
of [B, page 23].  We shall use such facts as $\bfone\forces``\dot P_{\eta,
\alpha}$ is a poset;'' see [S] for a proof.
By ``$\supt(p)$'' we mean $\{\beta\in\dom(p)\colon
p\restr\beta\notforces``p(\beta)=\bfone$''$\}$.
Notice that a stronger condition may have a support which is not
a superset of a weaker condition.

The following definition should be compared to [S 2, definition 1].

\proclaim Definition 1. Suppose that $\langle P_\eta\colon\eta<\alpha
\rangle$ is a forcing iteration. We define $Rlim^*(
P_\eta\colon\eta<\alpha)$. Let $\tilde P_\alpha$ be the
inverse limit of this sequence. We set $P_\alpha=Rlim^*(P_\eta\colon
\eta<\alpha)=\{p\in\tilde P_\alpha\colon(\forall q\leq_{\tilde P_\alpha}
p)(\forall\eta<\alpha)(\forall\dot\beta\in V^{P_\eta})
(\exists\zeta<\alpha)(\exists r\leq q\restr\zeta)(r\restr\eta=
q\restr\eta$ and $r\forces``${\rm if $\zeta<\dot\beta\leq\alpha$ then 
either $\cf(\dot\beta)
=\omega$ or $\supt(p\restr[\zeta,\dot\beta))$ is empty''$)\}$.}

If $P_\eta=Rlim^*(P_\xi\colon\xi<\eta)$ for every limit $\eta$, then
we say that $\langle
P_\eta\colon\eta\leq\alpha\rangle$ is an $RCS^*$-iteration.

Notice that $RCS^*$-iterations are forcing iterations in the
sense of [B], [S]. Thus we can use the results of [S, section 3]
to justify the usual abuses of notation which are used in
treatments of forcing iterations.

\proclaim Definition 2. Suppose $P$ is a poset  and
$\lambda$ is a sufficiently large regular cardinal and\/
$\mu$ is a cardinal (we allow $\mu=2$) and\/
$\{P,\mu\}\subseteq M\prec H_\lambda$ and $\vert M\vert=
\aleph_0$. 
Then we say that\/
$r$ is\/ {\rm $(M,\mu)$-pseudo-complete} iff whenever\/
$\sigma\in M^P$ then
there is $\gamma\in M\cap \mu$ such that\/
{\rm $r\forces``\sigma\in\check\mu$ {\rm implies} $\sigma=\check
\gamma$.''}

\proclaim Definition 3. $P$ is\/ {\rm $\mu$-pseudo-complete} iff
whenever\/ $\{P,\mu\}\subseteq M\prec H_\lambda$ and $\vert M\vert=
\aleph_0$ and $p\in P\cap M$, then there is $q\leq p$ such that
$q$ is $(M,\mu)$-pseudo-complete.

\proclaim Lemma 4. Suppose $P$ is $\mu$-pseudo-complete and
$G$ is a $V$-generic filter over $P$. Then $({}^\omega\mu)^V=
({}^\omega\mu)^{V[G]}$.

Proof: Immediate.

\proclaim Definition 5. Suppose $P$ is a $\mu$-pseudo-complete poset
and $\dot Q$ is a $P$-name for a poset. We say that\/ $\dot Q$
is\/ {\rm $\mu$-pseudo-complete relative to $P$} iff whenever
$\lambda$ is a sufficiently large regular cardinal and\/
$\{P*\dot Q,\mu\}\subseteq
 M\prec H_\lambda$ and $\vert M\vert=\aleph_0$ and $p=
(p_0,\dot p_1)\in P*\dot Q\cap M$ 
then there is $q=(q_0,\dot q_1)\leq p$ such that 
$q_0$ is $(M,\mu)$-pseudo-complete and
whenever
$\sigma\in M^{P*\dot Q}$ there is $\tau\in M^P$ such that\/ {\rm
$q_0\forces``\dot q_1\forces`\sigma\in\check\mu$ implies
$\sigma=\check\tau$'\thinspace'').}

Notice that the hypothesis on $\sigma$ and $\dot p_1$ is stronger than
requireing them to be in $M[G_P]$, and it is different from requiring
them to be in $\check M$.  Similarly the conclusion does not imply
that $\tau\in M[G_P]$ nor that $\tau\in\check M$.

\proclaim Definition 6. We say $P_\alpha$ is\/ {\rm strictly 
$\mu$-pseudo-complete} iff whenever\/ $\lambda$ is a sufficiently
large regular cardinal and\/ $\{P_\alpha,\mu\}\subseteq M\prec H_\lambda$
and\/ $\vert M\vert=\aleph_0$ and\/ $\eta\in\alpha\cap M$ and\/
$p\in P_\alpha\cap M$ and\/ $q\leq p\restr\eta$
and $q$ is $(M,\mu)$-pseudo-complete
then there is $r\leq p$
such that $r\restr\eta=q$ and\/
$r$ is $(M,\mu)$-pseudo-complete.

\proclaim Lemma 7. Suppose $P_\alpha$ is 
strictly $\mu$-pseudo-complete.
Then $P_\alpha$ is $\mu$-pseudo-complete.

Proof: Obvious (take $\eta=0$ in definition 6).

\proclaim Definition 8. We say that the
forcing iteration $\langle P_\eta\colon\eta\leq\kappa\rangle$
is\/ {\rm uniform} iff\/ 
{\rm $(\forall\alpha\leq\kappa)((\exists\gamma<\alpha)
(\bfone\forces_{P_\gamma}``\cf(\alpha)\leq\omega$'') or
$(\forall\gamma<\alpha)(\bfone\forces_{P_\gamma}``\cf(\alpha)>\omega$'')).}

The following is the main result of this paper.

\proclaim Theorem 9. Suppose $\langle P_\eta\colon\eta\leq\kappa\rangle$
is a uniform $RCS^*$-iteration and suppose $(\forall\eta<\kappa)
(\dot Q$ {\rm is $\mu$-pseudo-complete relative to $P_\eta)$.}
Then $P_\kappa$ is strictly $\mu$-pseudo-complete.

Proof: By induction on $\kappa$. Fix $M$, $\lambda$,
$p$, $\eta$, and $q$ as in definition 6.

Case 1: $\kappa$ is a successor ordinal.

Immediate.

Case 2: $\cf(\kappa) = \omega$.

Take $\langle\alpha_n\colon n\in\omega\rangle\in M$ an increasing
sequence cofinal in $\kappa$ with $\alpha_0=\eta$. Let $\langle\sigma_n\colon
n\in\omega\rangle$ enumerate $M^{P_\kappa}$. Build
$\langle r_n,\tau_n,q_n\colon n\in\omega\rangle$ by
recursion such that $q_0= q$ and $r_0=p$ and each of the
following holds:

($i$) $r_{n+1}\leq r_n$ and $r_{n+1}\restr\alpha_n= q_n$
and $r_{n+1}\restr[\alpha_n,\kappa)\in M^{P_{\alpha_n}}$ and
$r_{n+1}\forces``\sigma_n\in\check\mu$ implies $\sigma_n=\check
\tau_n$''

($ii$) $q_{n+1}\restr\alpha_n=q_n$ and $q_{n+1}\leq r_{n+1}\restr
\alpha_{n+1}$ and $q_n$ is $(M,\mu)$-pseudo-complete.

Now set $r = \bigcup\{q_n\colon n\in\omega\}$. Clearly this $r$ is
as required.

Case 3: $(\forall\gamma<\kappa)(\bfone\forces_{P_\gamma}``\cf(\kappa)>
\omega$'')

Build $\langle p_n,\tau_n,r_n,\zeta_n,q_n\colon
n\in\omega\rangle$ such that $p_0=p $ and $r_0=p$ and
$q_0=q$ and $\zeta_0=\eta$ and all of the following hold:

($i$) $\bfone\forces_{P_{\zeta_n}}``p_{n+1}\restr[\zeta_n,\kappa)
\forces`\sigma_n = \check\tau_{n+1}$'\thinspace'' and
$\tau_{n+1}\in M^{P_{\zeta_n}}$ and $p_{n+1}\leq r_n$ and $p_{n+1}\in M$
and $p_{n+1}\restr\zeta_n=r_n\restr\zeta_n$

($ii$) $r_{n+1}\restr\zeta_{n+1}
\forces``\supt(p_{n+1}\restr
[\zeta_{n+1},\kappa))=\emptyset$''
and $r_{n+1}\leq p_{n+1}$
and $r_{n+1}\restr\zeta_n=p_{n+1}\restr\zeta_n$ and $r_{n+1}\in M$
and $\zeta_{n+1}\in M$

($iii$) $q_{n+1}\leq r_{n+1}\restr\zeta_{n+1}$ and
$q_{n+1}\restr\zeta_n=q_n$ and $q_{n+1}$ is $(M,\mu)$-pseudo-complete

Now take $r\in P_\kappa$ such that $r\restr\zeta_n=q_n$ for every $n\in
\omega$ and $\supt(r)\subseteq\bigcup\{\zeta_n\colon n\in\omega\}$.
It is easy to see that this $r$ witnesses the desired conclusion.

Case 4: $\cf(\kappa)>\omega$ and $\bfone\forces_{P_\eta'}``cf(\kappa)=\omega$''

We may assume $\eta'\in M$. Furthermore, by increasing $\eta$ if necessary,
we may assume $\eta=\eta'$.
Let $X=\{\alpha<\kappa\colon\cf(\alpha)=\omega_1\}$.
Notice that $(\forall\gamma<\kappa)(\bfone\forces_{P_\gamma}``
\cf(\kappa)\geq\omega_2$''). This is because if there is some
$p'\in P_\gamma$ with $\gamma<\kappa$ and $p'\forces``\cf(\kappa)=
\omega_1$,'' then we may assume $\gamma\geq\eta$ and therefore
we have that $p'\forces``{}^\omega 2 \ne {}^\omega 2\cap V$,''
a contradiction.  Becuase $\cf(\kappa)$ is at
least $\omega_2$, we may take $\langle\dot\alpha_n
\colon n\in\omega\rangle\in M^{P_\eta}$ such that
$\bfone\forces``\{\dot\alpha_n\colon n\in\omega\}$ is a subset of 
$\check X$ which is cofinal in $\check\kappa$, and furthermore
$\check X\cap\dot\alpha_n$ is unbounded in $\dot\alpha_n$
for every $n\in\omega$.'' Notice that for every $p'\in
P_\eta$ and $\zeta$ such that $p'\forces``\dot\alpha_n=\check\zeta$''
we have that $(\forall\gamma<\zeta)(\forall q'\in P_\gamma)
(q'\restr\eta\leq p'$ implies $q'\forces``\cf(\zeta)\geq\omega_1$''),
because otherwise we would have an element of $P_\gamma$ which
forces that a real exists which is not in the ground model, which
we know to be impossible.

Now build $\langle p_n,\tau_n,r_n,\zeta_n,\dot\alpha'_n,
q_n\colon n\in\omega\rangle$ such that
$p_0=p$ and $r_0=p$ and $q_0=q$ and $\dot \alpha'_0=\zeta_0=\eta$
and all of the following hold:

($i$) $\bfone\forces_{P_{\zeta_n}}``p_{n+1}\restr[\zeta_n,\kappa)
\forces`\sigma_n=\check\tau_{n+1}$'\thinspace''
and $\tau_{n+1}\in M^{P_{\zeta_n}}$ and $p_{n+1}
\leq p_n$ and $p_{n+1}\in M$ and $p_{n+1}\restr\zeta_n
=r_n\restr\zeta_n$

($ii$) $r_{n+1}\restr\zeta_{n}
\forces``\supt(p_{n+1}\restr
[\zeta_{n},\zeta_{n+1})\forces`supt(p_{n+1}\restr[\zeta_{n+1},
\dot\alpha'_n)=\emptyset$'\thinspace''
and $\zeta_{n+1}\in M$  and $\cf(\zeta_{n+1})=\omega_1$
and $r_{n+1}\leq p_{n+1}$
and $r_{n+1}\restr\zeta_n=p_{n+1}\restr\zeta_n$ and $r_{n+1}\in M$

($iii$) $\bfone\forces_{P_{\zeta_n}}``\dot\alpha'_{n+1}={\rm max}
(\dot\alpha_{n+1},\check\zeta_{n+1})$'' and $\dot\alpha'_{n+1}\in
M^{P_{\zeta_n}}$

($iv$) $q_{n+1}\leq r_{n+1}\restr\zeta_{n+1}$ and
$q_{n+1}\restr\zeta_n=q_n$ and $q_{n+1}$ is $(M,\mu)$-pseudo-complete

Now let $\zeta=\bigcup\{\zeta_n\colon n\in\omega\}$. Define $r\in P_\kappa$
such that $r\restr\zeta_n=q_n$ for every $n\in\omega$ and
$r\restr\beta\forces``r(\beta)=p_n(\beta)$ if
${\rm max}(\zeta,\dot\alpha'_n\leq\beta<\dot\alpha'_{n+1}$.'\thinspace''
Then $r$ is as required.

The theorem is established.

We now investigate applications of theorem 9.
It is easy to see that if $\mu^*$ is a cardinal and
$P$ is $\mu$-pseudo-complete and
$\bfone\forces_P``\mu^*$ is measurable and $\dot Q$ 
is Prikry forcing on $\mu^*$'' then $\dot Q$ is 
$\mu$-pseudo-complete relative to $P$ for every
$\mu < \mu^*$.  It is also straightforward
to show that if CH holds and $P$ is $2$-pseudo-complete
and $\bfone\forces``\dot Q$ is Namba forcing,'' then
$\dot Q$ is 2-pseudo-complete relative to $P$ (just
follow the usual proof that Namba forcing adds no reals,
keeping track of which objects are in $M$; see, e.g., [J, pp.~289--291]).
Also, if CH holds and $P$ is $\omega$-pseudo-complete
and $\bfone\forces_P``\dot Q$
is the poset of closed countable subsets of a
stationary $S^*\subseteq S^2_0$,'' then
$\dot Q$ is $\omega$-pseudo-complete relative to $P$.

One may ask whether having two different versions of revised
countable support iteration is actually necessary.
It is not hard to show that [S 2, theorem 3] does not
hold for $RCS^*$, and hence it appears unlikely that the
$RCS^*$ iteration of semi-proper forcings is semi-proper.
Hence we cannot abandon $RCS$ in favor of $RCS^*$.
In the converse direction, Shelah proves that the $RCS$ iteration
of the three posets with which we are concerned do not
add reals, but only by invoking an extra assumption; roughly,
that certain $\sigma$-closed cardinal collapses are interspersed
throughout the course of the iteration.  It is problematic
whether this assumption can be removed under $RCS$.
Therefore, it may indeed be necessary to have
two versions of revised countable support,
depending on what sort of preservation
is required (especially preservation of semi-properness versus
preservation of not adding reals).

\medskip

To do:

(1) Exposition of the rest of [Sh, chapter XI] -- i.e.,
deriving solution to Friedman's problem, Avraham's problem,
and the precipitous filter problem from theorem 9

(2) Show that if $\langle P_\eta\colon\eta\leq\kappa\rangle$ is
an $RCS^*$ iteration (is uniformness necessary?) and $P_\kappa$
has $\lambda$-c.c.~and $\bfone\forces_{P_\kappa}``\vert A\vert
<\lambda$'' then $(\exists\alpha<\kappa)(A\in V^{P_\alpha})$.
Note that it suffices to establish that the ``Kunen axioms''
hold for $RCS^*$, similarly to the last section of [S 2].

(3) Show that if each $\dot Q_\eta$ is proper then $RCS^*$ is
``the same'' as countable support, similar to the penultimate
section of [S 2].

(4) Assuming (3) is correct, combine the ideas of [S, theorem 36]
and theorem 9 to obtain preservation of ``weakly 2-pseudo-complete plus
$<\omega_1$-proper'' under (uniform?) countable support iteration;
check whether this can replace ``${\cal D}$-completeness for some
2-complete completeness system ${\cal D}$'' in each of the
applications found in [Sh, chapter VIII],
using the ideas of [S, theorem 62] for each relevant poset.

\medskip

\medskip

\noindent{\bf References}

\medskip

[B] Baumgartner, J., ``Iterated Forcing'' in {\sl Surveys in Set Theory,}
A.R.D.~Mathias (ed.), Cambridge University Press, 1979

\medskip

[J] Jech, T., {\sl Set Theory,} Academic Press, 1978

\medskip

[S] Schlindwein, C.,
``Consistency of Suslin's hypothesis, a non-special Aronszajn tree,
and GCH,'' {\sl Journal of Symbolic Logic,} vol.~{\bf 59}, pp.~1--29, 1994

\medskip

[S 2] Schlindwein, C., ``Simplified RCS iterations,''
{\sl Archive for Mathematical Logic,} vol.~{\bf 32}, pp.~341--349, 1993

\medskip

[Sh] Shelah, S., {\sl Proper Forcing,} Lecture Notes in Mathematics
{\bf 940}, Springer-Verlag, 1982

\vfill\eject

\end{document}